\documentclass[oneside, a4paper,reqno]{amsart}
\usepackage{amscd,amssymb,amsmath,graphicx,verbatim,mathrsfs,xcolor}
\usepackage{wasysym}
\usepackage{hyperref}
\usepackage[TS1,OT1,T1]{fontenc}
\usepackage{lscape} 
\usepackage{fullpage} 
\usepackage[all]{xy}
\newtheorem{theorem}{Theorem}[section]
\newtheorem{lemma}[theorem]{Lemma}

\newtheorem{proposition}[theorem]{Proposition}

\newtheorem*{theorem*}{Theorem}

\theoremstyle{definition}
\newtheorem{definition}[theorem]{Definition}

\newtheorem{example}[theorem]{Example}

\theoremstyle{remark}
\newtheorem{remark}[theorem]{Remark}

\usepackage[pagewise]{lineno}

\newcommand{\Xcal}{\ensuremath{\mathcal{X}}}

\newcommand{\Tcal}{\ensuremath{\mathcal{T}}}
\newcommand{\Fcal}{\ensuremath{\mathcal{F}}}

\newcommand{\Acal}{\ensuremath{\mathcal{A}}}

\newcommand{\Ucal}{\ensuremath{\mathcal{U}}}
\newcommand{\Vcal}{\ensuremath{\mathcal{V}}}

\numberwithin{equation}{section}

\begin{document}
\title{Quantity vs.~size in representation theory}\
\author{Jorge Vit\'oria}
\address{Jorge Vit\'oria, Dipartimento di Matematica e Informatica, Universit\`a degli Studi di Cagliari, Palazzo delle Scienze, via Ospedale, 72, 09124 Cagliari, Italy}
\email{jorge.vitoria@unica.it}
\urladdr{https://sites.google.com/view/jorgevitoria/}
\begin{abstract}
In this note, we survey two instances in the representation theory of finite-dimensional algebras where the quantity of a type of structures is intimately related to the size of those same structures. More explicitly, we review the fact that (1) a finite-dimensional algebra admits only finitely many indecomposable modules up to isomorphism if and only if every indecomposable module is finite-dimensional; (2) the category of modules over a finite-dimensional algebra admits only finitely many torsion classes if and only if every torsion class is generated by a finite-dimensional module.
\end{abstract}
\thanks{The author is very grateful for many insightful discussions with Lidia Angeleri H\"ugel, Frederik Marks and Rosanna Laking on these topics, and for their comments on a preliminary version of the article. The author is also grateful to Jorge Freitas, Samuel Lopes and Diogo Oliveira e Silva for the invitation to write this survey.
}
\maketitle

\section{Introduction}
\label{sec:a}
In representation theory we strive to understand how the properties of a given ring are reflected on that ring's actions on abelian groups. These actions are formalised by the notion of a \textit{module over a ring}, and they are often difficult, if not impossible, to classify completely. Hence, questions in the subject area typically include the problem of classifying modules with certain common properties (simplicity, indecomposability, projectivity, injectivity, ...) up to isomorphism. One could call this a \textit{microscopic} approach to representation theory, in which the main actors are the actual modules over a given ring. These are difficult problems and, more often than not, their solutions involve a fair amount of combinatorics. Another strand of representation theory takes instead a bird's eye view of the subject, i.e. a \textit{macroscopic} point of view, considering the category of all modules and its subcategories as central objects of study. Typical questions within this line of thought  include classification problems for subcategories of modules subject to certain properties. Here the tools are more of a homological and categorical nature. 

These different points of view are used to study a broad range of rings, from group rings (in representation theory of finite groups) to universal enveloping algebras (in Lie theory) from commutative rings (often occurring in algebraic geometry) to path algebras of quivers (central to the study of finite-dimensional algebras). Recall that an (associative) algebra over a field $\mathbb{K}$ is nothing but a ring with a compatible $\mathbb{K}$-vector space structure. Here, we will be focused on some aspects of the representation theory of $\mathbb{K}$-algebras which are finite-dimensional as $\mathbb{K}$-vector spaces. Our aim is to survey two results that attempt to answer the following type of question:
\medskip

\noindent\textbf{Question:} To which extent do finite-dimensional modules over a finite-dimensional algebra $\Lambda$ control \textit{the structure} of the category of all $\Lambda$-modules? 

\medskip

To make this question precise, we need to establish what we mean by \textit{structure}. In this paper this expression will have two meanings. The first surveyed result is a classical theorem in representation theory, and it discusses when is it true that any given module can be built from finite-dimensional modules using the simplest operation available: direct sums. On the second result, however, \textit{the structures} we want to have under control are distinguished classes of modules, called torsion classes. We then aim to answer the question of whether any torsion class in the category of all modules is determined by the finite-dimensional modules contained in it.

In both settings, however, a remarkable pattern arises: the \textit{fewer} the objects under consideration in the finite-dimensional world (may they be modules or torsion classes), the tighter is the grip that finite-dimensional modules have on the whole category. In other words, \textit{quantity} controls \textit{size}. Moreover, and perhaps equally surprising, the converse also holds.

This note is structured as follows. In Section 2 we discuss some basic notions from representation theory of finite-dimensional algebras. In Section 3 we discuss, without proof, the first surveyed result: a well-known, classical theorem in representation theory, due to Auslander (\cite{Aus, Aus2}), Fuller-Reiten (\cite{FR}) and Ringel-Tachikawa (\cite{RT}) concerning algebras of finite representation type. We illustrate the properties under study through some examples. In Section 4 we look at torsion classes and discuss some examples. Finally, in Section 5 we survey a recent result from \cite{AMV4} that builds on the work of Demonet-Iyama-Jasso (\cite{DIJ}). It focuses on algebras whose module category admits only finitely many torsion classes. We sketch a proof of the theorem, leaving out some technical facts that we state without proof. In that respect, the choice made in this paper is to present the results from a torsion-theoretic point of view, leaving outside of the exposition the (intrinsic) relation of the arguments to $\tau$-tilting theory (\cite{AIR}) or silting modules~(\cite{AMV1}).

\medskip



\section{Representations of finite-dimensional algebras}

Throughout, let $\mathbb{K}$ be an \textbf{algebraically closed field} and $\Lambda$ a \textbf{finite-dimensional $\mathbb{K}$-algebra}. In this section we discuss finite-dimensional algebras and their representations. For a thorough introduction to the subject, we refer the reader to \cite{ASS} or \cite{P}. 

\begin{example}\label{A3}
Let $Q$ be a finite directed graph (usually called a \textbf{quiver}). Consider the $\mathbb{K}$-vector space spanned by all oriented paths in $Q$, and endow it with a multiplication defined by concatenation of paths when possible, and zero when concatenation is not possible. This yields a $\mathbb{K}$-algebra, denoted by $\mathbb{K}Q$, called the \textbf{path algebra} of $Q$. If $Q$ has no oriented cycles then $\mathbb{K}Q$ is finite-dimensional. For example, if $Q$ is the quiver
$$\xymatrix{1\ar[r]^\alpha&2\ar[r]^\beta&3}$$
then $\Lambda:=\mathbb{K}Q$ is a 6-dimensional $\mathbb{K}$-vector space spanned by the stationary paths $e_1$, $e_2$, $e_3$, and the paths $\alpha$, $\beta$ and $\beta\alpha$. Multiplication comes by concatenation as explained above: $\beta\cdot\alpha = \beta\alpha$, while $\alpha \cdot \beta=0$. This algebra  is in fact isomorphic to the algebra of $3\times 3$ lower triangular matrices over $\mathbb{K}$.
\end{example}

We will be looking at two categories associated to $\Lambda$. 
\begin{itemize}
\item $\mathsf{Mod}(\Lambda)$: the category of left $\Lambda$-modules;
\item $\mathsf{mod}(\Lambda)$: the category of finite-dimensional left $\Lambda$-modules.
\end{itemize}

If two finite-dimensional algebras have equivalent categories of modules, they are indistinguishable from a representation-theoretic standpoint. Algebras related in this way are said to be \textbf{Morita equivalent}. We will recall that any finite-dimensional algebra $\Lambda$ is Morita equivalent to a quotient of a path algebra. As a consequence, the category of modules over $\Lambda$ is equivalent to the category of bound representations of a quiver. 

A \textbf{representation of a quiver} Q over $\mathbb{K}$ is the assignment of a vector space to each vertex of $Q$ and a (compatibly chosen) linear map to each arrow of $Q$. A representation of $Q$ is \textbf{bound} by an ideal $I$ of $\mathbb{K}Q$ if the linear maps chosen for the arrows of the quiver compose and add up to the zero map when following a linear combination of paths contained in $I$. We denote by $\mathsf{Rep}(Q,I)$ the category of representations of $Q$ bound by $I$ (where morphisms between representations are given by linear maps at every vertex making the obvious diagrams commute). We refer to \cite[Corollary I.6.10 and Theorem II.3.7]{ASS} for the following theorem.

\begin{theorem}\label{rep}
Let $\Lambda$ be a finite-dimensional $\mathbb{K}$-algebra. Then there is a quiver $Q$ and an ideal $I$ of $\mathbb{K}Q$ such that $\mathsf{Mod}(\Lambda)\cong\mathsf{Mod}(\mathbb{K}Q/I)\cong \mathsf{Rep}(Q,I)$.
\end{theorem}

The equivalence above restricts to an equivalence between  finite-dimensional $\Lambda$-modules and finite-dimensional representations of $Q$ bound by $I$. In essence, the theorem indicates that the study of modules over any finite-dimensional $\mathbb{K}$-algebra is an upgrade of classical linear algebra over $\mathbb{K}$. 
\begin{remark}
The classical linear algebra theorem on Jordan normal forms can be regarded as a classification of the isomorphism classes of finite-dimensional representations over a quiver $Q$ with one vertex and one loop. Indeed, any square matrix $n\times n$ corresponds to an endomorphism of $\mathbb{K}^n$, and any conjugation by invertible matrices corresponds to  an isomorphism between the associated representations of $Q$. Moreover, Jordan blocks of a given matrix correspond to the indecomposable summands of the associated representation (see Section \ref{sec 3} for the notion of indecomposability).  Since $\mathbb{K}Q$ is isomorphic to the polynomial algebra $\mathbb{K}[X]$, much of classical linear algebra can be regarded as the study of finite-dimensional $\mathbb{K}[X]$-representations.
\end{remark}

\section{Finite representation type}\label{sec 3}

It is not very surprising that over a finite-dimensiona $\mathbb{K}$-algebra $\Lambda$, the structure of the category of finite-dimensional modules, $\mathsf{mod}(\Lambda)$, is much better understood than the structure of the category of all modules $\mathsf{Mod}(\Lambda)$. In fact, as stated in the following theorem (see for example \cite[Theorem I.4.10]{ASS}), objects in the category $\mathsf{mod}(\Lambda)$ can be completely described by a set of fundamental blocks: the finite-dimensional indecomposable modules. A $\Lambda$-module $M$ is said to be \textbf{indecomposable} if, whenever $M\cong M_1\oplus M_2$, then either $M_1$ or $M_2$ must be zero. In other words, a module is indecomposable if it admits no nontrivial direct summands.

\begin{theorem}[Krull-Remak-Schmidt]
Every finite-dimensional $\Lambda$-module $M$ is a direct sum of indecomposable $\Lambda$-modules, which are uniquely determined by $M$ up to isomorphism. 
\end{theorem}

The same statement, however, does not hold in general for infinite-dimensional $\Lambda$-modules. In fact, infinite-dimensional modules may exhibit a striking pathological property: not having \textbf{any} indecomposable direct summands! Such modules are called \textbf{superdecomposable} (see Example \ref{3K}). 

Let us consider the following four \textit{dream} properties for the representation theory of a finite-dimensional algebra $\Lambda$.
\begin{enumerate}
\item[(RF1)] Every $\Lambda$-module is a direct sum of indecomposable $\Lambda$-modules.
\item[(RF2)] Every indecomposable $\Lambda$-module is finite-dimensional.
\item[(RF3)] There are only finitely many indecomposable finite-dimensional $\Lambda$-modules up to isomorphism.
\item[(RF4)] There are only finitely many indecomposable $\Lambda$-modules up to isomorphism.
\end{enumerate}

The property (RF1) is called \textbf{pure semisimplicity} and the property (RF3) is called \textbf{representation-finiteness}. While (RF1) gives structural information on the category $\mathsf{Mod}(\Lambda)$, one can consider (RF2) a property regarding \textit{size} and both (RF3) and (RF4) properties regarding \textit{quantity}. We will establish a connection between all of these very soon.

\begin{example}\label{A3 again}
The algebra $\Lambda$ from Example \ref{A3} is representation-finite, i.e. $\Lambda$ satisfies (RF3). Up to isomorphism, there are precisely six indecomposable finite-dimensional $\Lambda$-modules and there are no infinite-dimensional indecomposables (thus $\Lambda$ satisfies also (RF2) and (RF4)). Using Theorem \ref{rep}, we can use quiver representations to describe these indecomposable modules:

$$P_3:=(\xymatrix{0\ar[r]&0\ar[r]&\mathbb{K}})\ \ \ \ \ S_2:=(\xymatrix{0\ar[r]&\mathbb{K}\ar[r]&0})\ \ \ \ \ \ S_1:=(\xymatrix{\mathbb{K}\ar[r]&0\ar[r]&0})$$
$$P_2:=(\xymatrix{0\ar[r]&\mathbb{K}\ar[r]^{1}&\mathbb{K}})\ \ \ \ \ P_1:=(\xymatrix{\mathbb{K}\ar[r]^{1}&\mathbb{K}\ar[r]^{1}&\mathbb{K}})\ \ \ \ \ I_2:=(\xymatrix{\mathbb{K}\ar[r]^{1}&\mathbb{K}\ar[r]&0}).$$

\noindent It can also be shown that every $\Lambda$-module is isomorphic to a direct sum of copies of these six $\Lambda$-modules, thus proving that $\Lambda$ also satisfies (RF1).
\end{example}

\begin{example}\label{3K}
It is very easy to produce examples of finite-dimensional algebras that do not satisfy (RF3) (and that, therefore, do not satisfy (RF4) either). For example, if $\Lambda$ is the path algebra over $\mathbb{K}$ of the quiver
$$Q\colon \xymatrix{1\ar@<1.7ex>@/^/[rr]^\alpha\ar@<-1.7ex>@/_/[rr]^\gamma\ar[rr]^\beta&&2}$$
we can produce infinitely many pairwise non-isomorphic indecomposable finite-dimensional $\Lambda$-modules. Given $\lambda$ in $\mathbb{K}$, the representation of $Q$
$$M_\lambda\colon \xymatrix{\mathbb{K}\ar@<1.7ex>@/^/[rr]^\lambda\ar@<-1.7ex>@/_/[rr]^1\ar[rr]^0&&\mathbb{K}}$$
is indecomposable. Since $\mathbb{K}$ is infinite (it is algebraically closed), we immediately get infinitely many indecomposable finite-dimensional $\Lambda$-modules. Moreover, it is easy to check that $M_\lambda$ is isomorphic to $M_\mu$ if and only if $\lambda=\mu$. We can also produce an indecomposable infinite-dimensional $\Lambda$-module, showing that $\Lambda$ does not satisfy (RF2). An explicit example is 
$$G\colon \xymatrix{\mathbb{K}(X)\ar@<1.7ex>@/^/[rr]^X\ar@<-1.7ex>@/_/[rr]^1\ar[rr]^0&&\mathbb{K}(X)}$$
where $\mathbb{K}(X)$ is the field of rational functions in one variable. This module belongs to an important family of indecomposable infinite-dimensional modules: the \textbf{generic modules}. These are indecomposable infinite-dimensional modules that, as modules over their own endomorphism ring, have finite length. They play an important role in controlling the overall representation theory of a finite-dimensional algebra. We refer to \cite{Ringel} and references therein for further information on generic modules.

Finally, we produce a superdecomposable $\Lambda$-module, thus showing (in a rather extreme way!) that (RF1) is not satisfied by $\Lambda$. This example can be found in \cite{Ringel}. We need two nontrivial ingredients: the existence of injective envelopes (see for example \cite[Corollary I.5.14]{ASS}) and the existence of a particularly nice functor between two categories of modules (see below).
\begin{itemize}
\item First consider the free algebra $\mathbb{K}\langle X,Y\rangle$ in two variables over $\mathbb{K}$. Note that this infinite-dimensional $\mathbb{K}$-algebra is isomorphic to the path algebra over $\mathbb{K}$ of the quiver with one vertex and two loops $X$ and $Y$ on that vertex. Let $I$ be the injective envelope of $\mathbb{K}\langle X,Y\rangle$ in $\mathsf{Mod}(\mathbb{K}\langle X,Y\rangle)$. We show that $I$ is superdecomposable. If $N\neq 0$ is a summand of $I$, then it intersects $\mathbb{K}\langle X,Y\rangle$  nontrivially  since $I$ is an essential extension of $\mathbb{K}\langle X,Y\rangle$. Let $a\neq0$ be an element in $N\cap \mathbb{K}\langle X,Y\rangle$ and consider the injective envelope $J$ of $\mathbb{K}\langle X,Y\rangle Xa$. It follows that $J$ is a nonzero summand of $N$ (by the injectivity of $J$). It remains to see that $J\neq N$, and this follows from the fact that the element $Ya$ of $N\cap \mathbb{K}\langle X,Y\rangle$ cannot lie in $J$ since $\mathbb{K}\langle X,Y\rangle Xa\cap \mathbb{K}\langle X,Y\rangle Ya=\{0\}$ and since $J$ is an essential extension of $\mathbb{K}\langle X,Y\rangle Xa$. Thus $N$ is not indecomposable.
\item There is a functor $F\colon \mathsf{Mod}(\mathbb{K}\langle X,Y\rangle)\longrightarrow \mathsf{Mod}(\Lambda)$ sending a (left) $\mathbb{K}\langle X,Y\rangle$-module $M$ to the representation of $Q$ given by 
$$F(M)\colon\xymatrix{M\ar@<1.7ex>@/^/[rr]^X\ar@<-1.7ex>@/_/[rr]^1\ar[rr]^Y&&M}$$
where $X$ and $Y$ are the linear maps resulting from the left action of $X$ and $Y$ on $M$. This functor is full and exact (\cite{Ringel}) and, as a consequence, $F(I)$ is a superdecomposable $\Lambda$-module. 
\end{itemize}
\end{example}

These examples suggest that the properties (RF1)-(RF4) come in a single package and cannot be satisfied separately. The following fundamental theorem in the representation theory of finite-dimensional algebras states that, indeed, these properties are equivalent. In other words, a finite-dimensional algebra has \textit{very few} (= finitely many) indecomposables up to isomorphism if and only if all indecomposables are \textit{small} (= finite-dimensional).

\begin{theorem}\cite{Aus, Aus2, FR, RT}
For a finite-dimensional $\mathbb{K}$-algebra $\Lambda$, the conditions (RF1), (RF2), (RF3) and (RF4) are equivalent.
\end{theorem}

\section{Torsion pairs}

Sometimes it is useful to have a birds' eye view of a category of modules and, rather than analysing the category module by module, organise collections of modules which share certain properties into certain classes. 

\begin{example}
Every abelian group has a subgroup given by the elements that have finite order. This is called the torsion subgroup. The quotient of an abelian group by its torsion subgroup yields an abelian group where no element has finite order. We may therefore say that the classes of torsion abelian groups and torsionfree abelian groups give us some valuable information on the category of abelian groups.
\end{example}

The following definition is an abstraction of the example above to an arbitrary abelian category.

\begin{definition}
A pair $(\Tcal,\Fcal)$ of full subcategories of an abelian category $\Acal$ is a \textbf{torsion pair} if
\begin{enumerate}
\item $\mathsf{Hom}_\Acal(T,F)=0$ for any $T$ in $\Tcal$ and any $F$ in $\Fcal$.
\item For any $X$ in $\Acal$, there are objects $t(X)$ and $f(X)$ in $\Tcal$ and $\Fcal$ respectively, and a short exact sequence of the form
$$0\longrightarrow t(X)\longrightarrow X\longrightarrow f(X)\longrightarrow 0$$
\end{enumerate}
Given a torsion pair $(\Tcal,\Fcal)$ in $\Acal$, we say that $\Tcal$ is a \textbf{torsion class} and $\Fcal$ a \textbf{torsionfree class}.
\end{definition}

It follows from the definition that the torsion class determines the torsionfree class. For a subcategory $\Xcal$ of $\mathsf{Mod}(\Lambda)$, denote by $\Xcal^\perp$ the full subcategory of $\mathsf{Mod}(\Lambda)$ whose objects are the modules $Y$ for which $\mathsf{Hom}_\Lambda(X,Y)=0$ for all $X$ in $\Xcal$. Dually, one may also define ${}^\perp\Xcal$. Given a torsion pair $(\Tcal,\Fcal)$, we have $\Fcal=\Tcal^\perp$ and $\Tcal={}^\perp\Fcal$. It can also be shown that in the category $\mathsf{Mod}(\Lambda)$, a full subcategory $\Xcal$ is a torsion class if and only if it is closed under coproducts (i.e. for any family of objects in $\Xcal$, its coproduct lies also in $\Xcal$), quotients (i.e. any quotient of a module in $\Xcal$ also lies in $\Xcal$) and extensions (i.e. in any short exact sequence with outer terms in $\Xcal$, the middle term must also belong to $\Xcal$). Dually, torsionfree classes are those full subcategories closed under products, submodules and extensions.

\begin{example}
Consider $\mathsf{Mod}(\mathbb{Z})$, the category of abelian groups. Let $\Tcal$ be the class of abelian groups for which every element has finite order and $\Fcal$ the class of abelian groups that have no elements of finite order. It is easy to check the axioms listed above showing that $(\Tcal,\Fcal)$ is a torsion pair. It can, furthermore be shown that $\Fcal$ is the subcategory of abelian groups which are subgroups of a product of a (possibly infinite) number of copies of $\mathbb{Q}$ and that $\Tcal={}^\perp \mathbb{Q}$. We then say that this torsion pair is \textbf{cogenerated} by $\mathbb{Q}$.
\end{example}
\begin{example}
Let $\Lambda$ be the algebra from Example \ref{A3}. Since every module in $\mathsf{Mod}(\Lambda)$ is isomorphic to a direct sum of the indecomposable $\Lambda$-modules listed in Example \ref{A3 again}, and since torsion classes are closed under coproducts and summands, every torsion class is determined by the set of indecomposable modules lying in it. Hence, every torsion class in $\mathsf{Mod}(\Lambda)$ is the closure under coproducts and their direct summands of some set of indecomposable $\Lambda$-modules $\Xcal$ that is closed under quotients and extensions. We denote this \textit{additive closure} by $\mathsf{Add}(\Xcal)$, and the complete list of torsion classes in $\mathsf{Mod}(\Lambda)$ (excluding $\{0\}$ and $\mathsf{Mod}(\Lambda)$) is
$$\Tcal_1:=\mathsf{Add}(\{S_2\})\ \ \Tcal_2:=\mathsf{Add}(\{P_3\})\ \ \Tcal_3:=\mathsf{Add}(\{S_1\})$$
$$\Tcal_4:=\mathsf{Add}(\{P_2,S_2\})\ \ \Tcal_5:=\mathsf{Add}(\{I_2,S_1\})\ \ \Tcal_6:=\mathsf{Add}(\{P_3,S_1\})$$
$$\Tcal_7:=\mathsf{Add}(\{P_3,P_2,S_2\})\ \ \Tcal_8:=\mathsf{Add}(\{S_2,I_2,S_1\})\ \ \Tcal_9:=\mathsf{Add}(\{P_1,I_2,S_1\})$$
$$\Tcal_{10}:=\mathsf{Add}(\{S_2,P_1,I_2,S_1\})\ \ \Tcal_{11}:=\mathsf{Add}(\{P_3,P_1,I_2,S_1\})$$
$$\Tcal_{12}:=\mathsf{Add}(\{P_2,P_1,S_2,I_2,S_1\})$$
We can order these classes by inclusion, obtaining the following Hasse quiver, where an arrow $A\rightarrow B$ denotes a strict inclusion $A\supsetneq B$ with no element $C$ such that $A\supsetneq C\supsetneq B$. The arrows in the Hasse diagram are linked to a process called \textbf{mutation} (see, for example, \cite{AIR}). While we will not discuss this process, it plays an important role in the proof of Theorem~\ref{main thm}. 

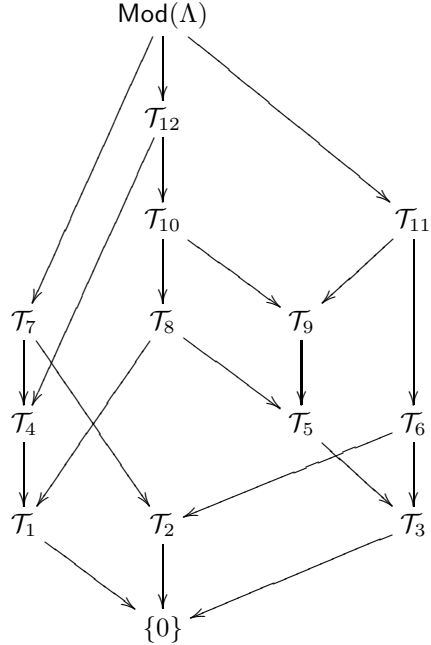
$$\xymatrix{&\mathsf{Mod}(\Lambda)\ar[lddd]\ar[d]\ar[rrdd]\\ &\Tcal_{12}\ar[dddl]\ar[d]\\ &\Tcal_{10}\ar[d]\ar[dr]&&\Tcal_{11}\ar[dl]\ar[dd]\\ \Tcal_7\ar[ddr]\ar[d] &\Tcal_8\ar[dr]\ar[ddl]&\Tcal_9\ar[d]\\ \Tcal_4\ar[d]&&\Tcal_5\ar[dr]&\Tcal_6\ar[d]\ar[dll]\\ \Tcal_1\ar[dr]&\Tcal_2\ar[d]&&\Tcal_3\ar[dll]\\ &\{0\}}$$

\noindent Note also that this Hasse quiver depicts a well-known object: the three-dimensional associahedron. For further details on the combinatorics of torsion pairs of path algebras of quivers we refer to \cite{IT}.

\end{example}
\begin{example}
For a $\Lambda$-module $M$, the subcategory $M^\perp$ of $\mathsf{Mod}(\Lambda)$ is closed under products, submodules and extensions and, thus, $M^\perp$ is a torsionfree class. The corresponding torsion class is then necessarily given by ${}^\perp(M^\perp)$, and it is clear that $M$ lies in ${}^\perp(M^\perp)$. In general it is not easy to describe which modules lie in ${}^\perp(M^\perp)$.  However, if $M$ is, for example, a projective $\Lambda$-module, then one can show that ${}^\perp(M^\perp)$ coincides with the subcategory  $\mathsf{Gen}(M)$ formed by all $\Lambda$-modules which are quotients of some coproduct of copies of $M$. In Section \ref{sec final} we describe the torsion pairs that are of the form $(\mathsf{Gen}(M),M^\perp)$ for a finite-dimensional $\Lambda$-module $M$.
\end{example}

\section{Torsion-finiteness}\label{sec final}

In this section we will discuss categories of modules that have only finitely many torsion classes. Let us first look at torsion classes of $\mathsf{Mod}(\Lambda)$, for a finite-dimensional $\mathbb{K}$-algebra $\Lambda$, which are of the form $\mathsf{Gen}(M)$ for some $\Lambda$-module $M$. They satisfy the following useful property.

\begin{lemma}\label{union}\cite[Lemma 3.10]{DIJ}
Suppose that $\Tcal=\mathsf{Gen}(M)$ is a torsion class. Suppose that there is an ascending sequence of torsion classes 
$$\Tcal_1\subseteq \Tcal_2\subseteq \Tcal_3\subseteq \cdots\subseteq \Tcal_n\subseteq \Tcal_{n+1}\subseteq \cdots$$
such that $\bigcup\limits_{i\geq 1}\Tcal_i=\Tcal$. Then the sequence stabilises.
\end{lemma}
\begin{proof}[Proof]
Given $n\geq 1$ such that $M$ lies in the torsion class $\Tcal_n$, all coproducts of copies of $M$ and their quotients must lie in $\Tcal_n$, proving that $\Tcal_n=\Tcal$.
\end{proof}

In the category $\mathsf{mod}(\Lambda)$ of finite-dimensional modules, quotients of finite coproducts of a finite-dimensional module $M$ sometimes also form a torsion class. Such a subcategory is denoted by $\mathsf{gen}(M)$. Clearly such subcategories also satisfy the property of Lemma \ref{union}. Torsion pairs in $\mathsf{mod}(\Lambda)$ and torsion pairs in $\mathsf{Mod}(\Lambda)$ are related by the following theorem. Given a subcategory $\Xcal$ of $\mathsf{mod}(\Lambda)$, denote by $\varinjlim \Xcal$ the subcategory of $\mathsf{Mod}(\Lambda)$ whose objects are direct limits of direct systems with terms in $\Xcal$. Note that, since a direct limit of a direct system is a quotient of the coproduct of the terms in that system, torsion classes in $\mathsf{Mod}(\Lambda)$ are closed under direct limits.

\begin{theorem}\label{CB}\cite[Lemma 4.4]{CB}
If $(\Ucal,\Vcal)$ is a torsion pair in $\mathsf{mod}(\Lambda)$, then $(\varinjlim \Ucal,\varinjlim \Vcal)$ is a torsion pair in $\mathsf{Mod}(\Lambda)$, and 
$(\varinjlim \Ucal,\varinjlim \Vcal)=(\mathsf{Gen}(\Ucal),\Ucal^\perp).$
This assignment is injective since
$$(\Ucal,\Vcal)=((\varinjlim\Ucal)\cap \mathsf{mod}(\Lambda), (\varinjlim\Vcal)\cap \mathsf{mod}(\Lambda)),$$
and a torsion pair $(\Tcal,\Fcal)$ in $\mathsf{Mod}(\Lambda)$ arises in this way if and only if $\Fcal$ is closed under direct limits in $\mathsf{Mod}(\Lambda)$.
\end{theorem}

The theorem above establishes a close relation between torsion classes of the form $\mathsf{gen}(M)$ and $\mathsf{Gen}(M)$ in $\mathsf{mod}(\Lambda)$ and in $\mathsf{Mod}(\Lambda)$ respectively, for a finite-dimensional $\Lambda$-module $M$. In fact,  it follows that $\mathsf{gen}(M)$ is a torsion class in $\mathsf{mod}(\Lambda)$ if and only if $\mathsf{Gen}(M)$ is a torsion class in $\mathsf{Mod}(\Lambda)$, in which case $\mathsf{Gen}(M)=\varinjlim \mathsf{gen}(M)$. The following result characterises the torsion pairs that are of the form $(\mathsf{Gen}(M), M^\perp)$. Recall that a \textbf{pure submodule} $Y$ of a $\Lambda$-module $X$ is a submodule such that for any right $\Lambda$-module $Z$, we have that $Z\otimes_\Lambda Y$ is still a submodule of $Z\otimes_\Lambda X$. For example, it is easy to check that if $X/Y$ is a flat module, then $Y$ is a pure submodule of $X$.

\begin{proposition}\label{new tp gen}
Let $\Lambda$ be a finite dimensional $\mathbb{K}$-algebra and let $(\Tcal,\Fcal)$ be a torsion pair in $\mathsf{Mod}(\Lambda)$. The following statements are equivalent.
\begin{enumerate}
\item $\Tcal\cap \mathsf{mod}(\Lambda)=\mathsf{gen}(M)$ for a finite dimensional $\Lambda$-module $M$;
\item $(\Tcal,\Fcal)=(\mathsf{Gen}(M),M^\perp)$ for a finite dimensional $\Lambda$-module $M$;
\item $\Tcal$ is closed under products and $\Fcal$ is closed under direct limits.
\end{enumerate}
\end{proposition}
\begin{proof}[Proof] Let $(\Tcal,\Fcal)$ be a torsion pair in $\mathsf{Mod}(\Lambda)$.

\medskip

(1) $\Rightarrow$ (2): Suppose that $\Tcal\cap\mathsf{mod}(\Lambda)=\mathsf{gen}(M)$. Since $\Tcal\cap\mathsf{mod}(\Lambda)$ is a torsion class in $\mathsf{mod}(\Lambda)$ it follows that, as seen above, $\mathsf{Gen}(M)=\varinjlim (\Tcal\cap\mathsf{mod}(\Lambda))$ is a torsion class, and that $\mathsf{Gen}(M)\subseteq \Tcal$. To prove that $\Tcal\subseteq \mathsf{Gen}(M)$ we need the following facts.
\begin{enumerate}
\item[(P1)] \cite[2.2, Example 3]{CB2} Every $\Lambda$-module is a pure submodule of the product of its finite-dimensional quotients.
\item[(P2)] \cite[Theorem 4.2]{CB} For a torsion class $\Ucal$ in $\mathsf{mod}(\Lambda)$, $\varinjlim \Ucal$ is always closed under pure submodules and, moreover, it is  closed under products if and only if $\Ucal=\mathsf{gen}(N)$ for a finite-dimensional $\Lambda$-module $N$.
\end{enumerate}
If $X$ is a module in $\Tcal$ then, by (P1), $X$ is a pure submodule of its finite-dimensional quotients, all of which lie in $\Tcal\cap\mathsf{mod}(\Lambda)=\mathsf{gen}(M)$. Since, by (P2), $\mathsf{Gen}(M)=\varinjlim \mathsf{gen}(M)$ is closed under products and pure submodules, $X$ lies in $\mathsf{Gen}(M)$.

\medskip

(2) $\Rightarrow$ (3): Since $M$ is finite-dimensional,  the functor $\mathsf{Hom}_\Lambda(M,-)$ commutes with direct limits and, thus, $\Fcal=M^\perp$ is closed under direct limits. Moreover, since $\mathsf{Gen}(M)\cap \mathsf{mod}(\Lambda)=\mathsf{gen}(M)$, it follows from the statement (P2) that $\mathsf{Gen}(M)$ is closed under products.

\medskip

(3) $\Rightarrow$ (1): If $\Fcal$ is closed under direct limits then, by Theorem \ref{CB}, 
$$(\Tcal,\Fcal)=(\varinjlim (\Tcal\cap \mathsf{mod}(\Lambda)), \varinjlim (\Fcal\cap \mathsf{mod}(\Lambda))).$$
By the fact (P2) cited above, it then follows that there is a finite-dimensional $\Lambda$-module $M$ such that $\Tcal\cap \mathsf{mod}(\Lambda)=\mathsf{gen}(M)$.
\end{proof}

If one wishes to classify all torsion pairs in $\mathsf{Mod}(\Lambda)$, the following four \textit{dream} properties could be of help, just like in Section \ref{sec 3}.

\begin{enumerate}
\item[(TF1)] Every torsion pair in $\mathsf{Mod}(\Lambda)$ is of the form $(\varinjlim \Ucal,\varinjlim \Vcal)$ for a torsion pair $(\Ucal,\Vcal)$ in $\mathsf{mod}(\Lambda)$;
\item[(TF2)]  Every torsion class in $\mathsf{Mod}(\Lambda)$ is of the form $\mathsf{Gen}(M)$  for a finite-dimensional $\Lambda$-module $M$;
\item[(TF3)] There are only finitely many torsion classes in $\mathsf{mod}(\Lambda)$;
\item[(TF4)] There are only finitely many torsion classes in $\mathsf{Mod}(\Lambda)$.
\end{enumerate}

Note that, once again, (TF1) is a property concerning the structure of torsion pairs in the category of $\Lambda$-modules, while (TF2) can be seen as a measure of \textit{size} and (TF3) and (TF4) as measures of \textit{quantity}. Just like in the previous section, it turns out that these properties (TF1)--(TF4) are equivalent to each other. In other words, a finite-dimensional algebra admits \textit{very few} (= finitely many) torsion classes in its category of (finite-dimensional) modules if and only if all torsion classes in its module category are generated by \textit{small} (= finite-dimensional) modules. The following theorem is essentially proved in \cite{AMV4}, using the language of support $\tau$-tilting modules (\cite{AIR}) and silting modules (\cite{AMV1}). The proof we present here purposefully avoids mentioning these modules, re-working the existing arguments from a torsion-theoretic point of view.

\begin{theorem}\cite[Theorem 4.8]{AMV4}\label{main thm}
For a finite-dimensional $\mathbb{K}$-algebra $\Lambda$, the conditions (TF1), (TF2), (TF3) and (TF4) are equivalent.
\end{theorem}
\begin{proof}[Proof]
In order to prove this theorem we need the following additional fact.
\begin{itemize}
\item[(P3)] \cite[Proposition 3.8]{DIJ} There are only finitely many torsion classes in $\mathsf{mod}(\Lambda)$ if and only if every torsion class in $\mathsf{mod}(\Lambda)$ is of the form $\mathsf{gen}(M)$ for a finite-dimensional $\Lambda$-module $M$.
\end{itemize}
This fact partially relies on the combinatorial technique mentioned in the previous section: mutation. Indeed, if there are either infinitely many torsion classes in $\mathsf{mod}(\Lambda)$ or if there is a torsion class which is not of the form $\mathsf{gen}(M)$, one can build a non-stabilising ascending chain of torsion classes 
$$\Tcal_1\subsetneq \Tcal_2\subsetneq \Tcal_3\subsetneq \cdots\subsetneq \Tcal_n\subsetneq \Tcal_{n+1}\subsetneq \cdots$$
which, with the help of Lemma \ref{union}, guarantees the other condition. We refer the reader to \cite{AIR} and \cite{DIJ} for further details.

\medskip

(TF1) $\Rightarrow$ (TF2): Let $(\Tcal,\Fcal)$ be a torsion pair in $\mathsf{Mod}(\Lambda)$. By assumption, $\Tcal=\varinjlim \Ucal$  and $\Fcal=\varinjlim \Vcal$, where $\Ucal=\Tcal\cap\mathsf{mod}(\Lambda)$ and $\Vcal=\Fcal\cap\mathsf{mod}(\Lambda)$. Let us first show that $\Tcal={}^\perp\Vcal$. Indeed, it is easy to check that ${}^\perp\Vcal$ is a torsion class containing $\Ucal$ and that the corresponding torsionfree class $({}^\perp\Vcal)^{\perp}$ contains $\Vcal$. Since every torsion class is closed under direct limits we immediately conclude that $\Tcal=\varinjlim\Ucal\subseteq {}^\perp\Vcal$. Since, by assumption, $\Fcal$ is also closed under direct limits, we have that $\Fcal=\varinjlim\Vcal\subseteq ({}^\perp\Vcal)^\perp$ and, therefore, we can also conclude that $\Tcal={}^\perp\Fcal\supseteq {}^\perp(({}^\perp\Vcal)^\perp)={}^\perp\Vcal$. This proves that $\Tcal={}^\perp \Vcal$, as wanted.

Now, since $\Vcal$ is a subcategory of finite-dimensional $\Lambda$-modules, ${}^\perp\Vcal$ is closed under products (see, for example, \cite[Example 2.3]{SD}). Therefore, (TF2) follows from Proposition \ref{new tp gen}. 

\medskip

(TF2) $\Rightarrow$ (TF3): If $\Ucal$ is a torsion class in $\mathsf{mod}(\Lambda)$, then by Theorem \ref{CB}, $\varinjlim \Ucal$ is a torsion class in $\mathsf{Mod}(\Lambda)$. From (TF2) we conclude that  $\varinjlim \Ucal=\mathsf{Gen}(M)$ for a finite-dimensional $\Lambda$-module $M$ and, thus, $\Ucal=\mathsf{Gen}(M)\cap\mathsf{mod}(\Lambda)=\mathsf{gen}(M)$. Finally, (TF3) follows (P3).

\medskip

(TF3) $\Rightarrow$ (TF1): This follows from (P3) and Proposition \ref{new tp gen}.

\medskip

(TF4) $\Rightarrow$ (TF3): This is a direct consequence of Theorem \ref{CB}.

\medskip

(TF2) $\Rightarrow$ (TF4): (TF2) implies that the assignment in Theorem \ref{CB} is a bijection (since $M^\perp$ is closed under direct limits for any finite-dimensional $\Lambda$-module $M$), i.e. there are as many torsion classes in $\mathsf{mod}(\Lambda)$ as in $\mathsf{Mod}(\Lambda)$. Since (TF2) is proved to be equivalent to (TF3), we conclude that there are finitely many torsion classes in $\mathsf{Mod}(\Lambda)$.
\end{proof}

Note that the properties (P1), (P2) and (P3) are fundamental to our proof, and they depend heavily on the fact that we are working with finite-dimensional algebras. We should, therefore, be very careful with any attempts to naively generalise the result above to larger classes of rings. We finish this survey with an example of a commutative noetherian (but not artinian) ring where the philosophy of this last section fails.

\begin{example}
Let $R$ be the (commutative, noetherian) ring of fractions of $\mathbb{Z}$ obtained by inverting all odd integers, i.e. 
$$R=\{\frac{a}{b}\in\mathbb{Q}\colon \mathsf{gcd}(b,2)=1\}.$$
This ring is a principal ideal domain, and its nontrivial ideals are just the powers of the maximal ideal $\mathfrak{p}$ generated by $2$. Finitely generated modules over a principal ideal domain are very well-understood. In our case, for every finitely generated $R$-module there is an isomorphism 
$$M\cong\bigoplus_{k\geq 0}(R/\mathfrak{p}^k)^{n_k(M)},$$ 
where $n_k(M)\neq 0$ for only finitely many $k$ (and $\mathfrak{p}^0=\{0\})$.
\begin{itemize}
\item We first show that $\mathsf{mod}(R)$ has only two nontrivial torsion classes. Let $\Tcal\neq \{0\}$ be a torsion class in $\mathsf{mod}(R)$, the subcategory of finitely generated $R$-modules. If there is a module $M$ in $\Tcal$ such that $M$ is faithful (i.e., such that $n_0(M)\neq 0$), then since $\Tcal$ is closed under direct summands, $R$ lies in $\Tcal$ and $\Tcal=\mathsf{mod}(R)$. If, on the other hand, every module in $\Tcal$ has a nonzero annihilator, it is easy to show that $R/\mathfrak{p}$ lies in $\Tcal$ (since it is a quotient of any nonzero module). Finally, observe that since $R/\mathfrak{p}^n$ is an iterated extension of $R/\mathfrak{p}$, we get that 
$$\Tcal=\{M\in\mathsf{mod}(R)\colon n_0(M)=0\}.$$ 

\item We now produce a torsion class in $\mathsf{Mod}(R)$ that is not generated by a finitely generated $R$-module. Since $R$ is a principal ideal domain, an $R$-module is injective if and only if it is divisible. It can then be shown that the injective $R$-modules are those in $\mathsf{Gen}(\mathbb{Q}\oplus \mathbb{Q}/R)$, and that they form a torsion class. Note additionally that $\mathbb{Q}$ is not a finitely generated $R$-module. Finally, it can be shown that any other generator of this same torsion class must contain $\mathbb{Q}$ as a summand.
\end{itemize} 
In conclusion, $\mathsf{mod}(R)$ satisfies the analogous condition to (TF3) and, yet, $\mathsf{Mod}(R)$ admits a torsion class that cannot be generated by a finitely generated $R$-module, therefore not satisfying the analogous condition to (TF2).
\end{example}

\end{document}